\documentclass[a4paper, reqno]{amsart}
\usepackage[utf8]{inputenc}
\usepackage[T2A]{fontenc}
\usepackage{amsmath,amssymb,amsthm}
\usepackage[a4paper,hmargin=2.5cm,vmargin=2.5cm]{geometry}
\usepackage[english]{babel}
\usepackage[unicode=true]{hyperref}
\usepackage{amscd}
\usepackage[all]{xy}
\usepackage{yfonts}
\usepackage{enumitem}
\usepackage{extarrows}

\usepackage{amsfonts}
\usepackage{graphicx}
\usepackage{mathrsfs}
\usepackage{amsmath}
\usepackage{amssymb}
\usepackage{pstricks}
\usepackage{verbatim}
\usepackage{mathtools}

\def\colon{\mathpunct{:}}

\title{Highly divisible cycles in homological Atiyah---Hirzebruch spectral sequences}

\author{Vasilii Rozhdestvenskii}
\thanks{This work was supported by the HSE University Basic Research Program and by the Theoretical Physics and Mathematics Advancement Foundation ``BASIS'' (grant no.\ 22-7-2-10-4). 
}
\address{{\sloppy
\parbox{0.9\textwidth}{
Steklov Mathematical Institute of Russian Academy of Sciences, Moscow, Russian Federation
\\[5pt]
National Research University Higher School of Economics, Russian Federation
}\smallskip}}
\email{vrozhd@mi-ras.ru  \medskip}

\begin{document}

\newtheorem{Th}{Theorem}
	\newtheorem*{Th*}{Theorem} 
	\newtheorem{Lem}{Lemma}
	\newtheorem*{Lem*}{Lemma}
	\newtheorem{Cor}{Corollary}
	\newtheorem*{Cor*}{Corollary}
	\newtheorem*{Frame}{Frame}
	\newtheorem*{Denom}{Denomination}
	\newtheorem*{Prop*}{Proposition}
	\newtheorem{Prop}{Proposition}
	\theoremstyle{definition}
	\newtheorem{definition}{Definition}
	\newtheorem*{definition*}{Definition}
	\newtheorem{Rem}{Remark}
	\newtheorem*{Rem*}{Remark}
	\newtheorem{Conv}{Convention}

 \renewcommand{\thesubsection}{\arabic{subsection}}

	\makeatletter
\newcommand{\colim@}[2]{%
  \vtop{\m@th\ialign{##\cr
    \hfil$#1\operator@font colim$\hfil\cr
    \noalign{\nointerlineskip\kern1.5\ex@}#2\cr
    \noalign{\nointerlineskip\kern-\ex@}\cr}}%
}
\newcommand{\colim}{%
  \mathop{\mathpalette\colim@{\rightarrowfill@\textstyle}}\nmlimits@
}
\makeatother
	
	\def\rk{\mathop{\mathrm{rk}}}
	\def\id{\mathord{\mathrm{id}}}
	\def\pt{\mathord{\mathrm{pt}}}
	\def\Ker{\mathop{\mathrm{Ker}}}
	\def\Im{\mathop{\mathrm{Im}}}
	\def\Tor{\mathop{\mathrm{Tor}}}
	\def\Ext{\mathop{\mathrm{Ext}}}
	\def\Tors{\mathop{\mathrm{Tors}}}
	\def\ord{\mathop{\mathrm{ord}}}
	\def\ex{\mathop{\mathrm{ex}}}
	\def\Indet{\mathop{\mathrm{Indet}}}
	\def\Hom{\operatorname{Hom}}
	\def\Ann{\operatorname{Ann}}
	\def\MU{\mathord{\mathrm{MU}}}
	\def\bu{\mathord{\mathrm{bu}}}
	\def\MSO{\mathord{\mathrm{MSO}}}
	\def\SO{\mathord{\mathrm{SO}}}
	\def\A{\mathord{\mathcal{A}}}
	\def\Sq{\mathord{\mathrm{Sq}}}
	\def\U{\mathord{\mathrm{U}}}
	\def\MG{\mathord{\mathrm{MG}}}
	\def\ch{\mathord{\mathrm{ch}}}
	\def\BU{\mathord{\mathrm{KU}}}
	\def\BP{\mathord{\mathrm{BP}}}
	\def\St{\mathord{\mathrm{St}}}
	\def\K{\mathord{\mathrm{K}}}
	\def\d{\mathord{\delta}}
	\newcommand{\Z}{\mathbb{Z}}
	\newcommand{\R}{\mathbb{R}}
	\newcommand{\Q}{\mathbb{Q}}
	\newcommand{\F}{\widehat{F}}
	\newcommand{\E}{\widehat{E}}
	\newcommand{\B}{\widehat{B}}
	\newcommand{\Zz}{\widehat{Z}}
	\newcommand{\dd}{\widehat{d}}
	\let\le\leqslant
	
	\let\ge\geqslant

\begin{abstract}
Let $F_*$ be a homology theory of finite type and let $\{E^r_{s,t}(X)=\frac{Z^r_{s,t}(X)}{B^r_{s,t}(X)}; d^r_{s,t}\}$ denote the Atiyah---Hirzebruch spectral sequence for $F_*$ of a $CW$-complex $X$. In this paper we  study elements of the maximal possible order in the groups $E^2_{n,0}(X)/Z^{r}_{n,0}(X)$. 
\end{abstract}
\maketitle

The initial motivation for the questions discussed below lies in the Steenrod problem on realization of cycles. Recall that in the late 1940s N. Steenrod posed  the following question (see \cite[Prob. 25]{Steenrod}): \textit{
Given a homology class $x\in H_n(X;\Z)$ of a topological space $X$, does there exist a closed oriented manifold~$M$ and \linebreak a continuous map $f\colon M \to X$ such that $f_*[M]=x$?} Here we consider this question under the assumption that $M$ is smooth. Fundamental results obtained by R. Thom (see \cite[Ch.~3]{Thom})  state that there exist Steenrod-non-realizable classes and an arbitrary integral homology class can be realized after multiplication by some positive integer, which can be chosen to depend only on its dimension. In connection with these results the following problem arose: 
\textit{
For each dimension~$n$ find the smallest positive integer $k_{\SO}(n)$ such that every integral $n$-dimensional homology class can be realized in the sense of Steenrod after multiplication by~$k_{\SO}(n)$.
}

The notion of Steenrod-realizability can be reformulated in terms of comparing homology theories. So, let $\MSO_*$ denote the oriented bordisms theory and let $\mu\colon \MSO_{n}(X) \to H_n(X;\Z)$ denote the canonical natural transformation. Then a class $x\in H_n(X;\Z)$ is Steenrod-realizable iff it lies in the image of $\mu$. Such a definition also motivates us to introduce a number $k_{\SO}^s(n)$ which is the smallest positive integer such that every homology class $x\in H_n(X; \Z)$ of a \textit{connective spectrum} $X$ becomes Steenrod-realizable after multiplication by $k_{\SO}^s(n)$. These numbers were independently computed by V.~M. Buchstaber in~\cite{Buch.1}, \cite{Buch.2} and G.~Brumfiel in \cite{Brumfiel} and they obtained that 
$
k^s_{\SO}(n)=\prod\limits_{p\ \text{prime},\ p>2}p^{\left[\frac{n-1}{2(p-1)}\right]}.
$
The standard comparison of stable and unstable categories implies that $k_{\SO}(n)$ divides $k^{s}_{\SO}(n-1)$ and $k^{s}_{\SO}([\frac{n-2}{2}])$ divides $k_{\SO}(n)$. So the Brumfiel---Buchstaber result also provides (multiplicative) upper and lower bounds on $k_{\SO}(n)$. \linebreak  A significant improvement of the Brumfiel---Buchstaber lower bound was obtained in the author's work~\cite{R}, where it was proved that $k_{\SO}(n)$ is divisible by the maximal odd divisor of $\left[\frac{n-1}{2}\right]\hskip-1mm \mathord{\vcenter{\hbox{\LARGE\textup !}}}$ and this lower bound is asymptotically equivalent in the logarithmic scale to the Brumfiel---Buchstaber upper bound (see \cite[Th.1 and Cor.2]{R}). Moreover, for $n<24$ the number $k_{\SO}(n)$ is equal to the author's lower bound (\cite[Th.3]{R}). 

In this paper we are interested in analogous questions for an arbitrary homology theory. Namely, let~$F$ be a spectrum, $\widetilde{F}$ be the connective cover of $F$, and $\mu_{F}\colon \widetilde{F}~\to~H\pi_0(F)$ be the canonical map to the Eilenberg---MacLane spectrum $H\pi_0(F)$. Then $\mu_{F}$ induces a natural transformation of the corresponding homology theories, which we will also denote by $\mu_{F}\colon \widetilde{F}_n(X) \to H_n(X;\pi_{0}F)$. Now one can define by analogy the numbers $k_{F}(n)$ and $k_{F}^s(n)$ to be the smallest positive integers such that every homology class $x\in H_n(X; \pi_{0}F)$ of a $CW$-complex (respectively, connective spectrum) $X$ lies in the image of $\mu_{F}$ after multiplication by $k_{F}(n)$ (respectively,  by $k_{F}^s(n)$). Note that in the case $F=\MSO$ we have that $k_{\SO}(n)=k_{F}(n)$ and $k_{\SO}^s(n)=k_{F}^s(n)$.  

For a general spectrum $F$ the problem of computing both the numbers $k_{F}(n)$ and $k^s_{F}(n)$ is highly non-trivial. Nevertheless, whereas the numbers $k^s_{F}(n)$ can in certain cases be computed (for example, \linebreak when $F=\MSO, \MU$ or $\BU$), the complete determination of the numbers $k_{F}(n)$
remains, in general, an intractable problem with no visible approaches (see Section $4$ for some reasons for this). In view of this, we focus on obtaining as much information as possible about the numbers $k_{F}(n)$ under the assumption that we know the numbers $k^{s}_{F}(n)$. To make this precise, let $k_{F}(n,l)$ denote the minimal positive integer such that for any $n$-dimensional homology class $x\in H_n\bigl(X; \pi_{0}F\bigr)$ of an  $l$-connected $CW$-complex~$X$, the class $k_{F}(n,l)x$ lies in the image of $\mu_{F}$. Observe that the introduced numbers are related in the following way: 
\begin{gather}
\begin{aligned}
\xymatrix{
k^{s}_{F}(n-1) \ar@{->>}[d] \ar@{->>}[r] & k^{s}_{F}(n-2) \ar@{->>}[d] \ar@{->>}[r] & \cdots \ar@{->>}[d] \ar@{->>}[r] & k^{s}_{F}([\frac{n-2}{2}]) \ar@{=}[d] \ar@{->>}[r] & \cdots \ar@{=}[d]  \ar@{->>}[r] & k^{s}_{F}(0)=1\ar@{=}[d]
\\
k_{F}(n)=k_{F}(n,0) \ar@{->>}[r] & k_{F}(n,1) \ar@{->>}[r] & \cdots \ar@{->>}[r] & k_{F}(n, [\frac{n+1}{2}])  \ar@{->>}[r] & \cdots \ar@{->>}[r] & k_{F}(n,n-1)=1, 
}
\end{aligned}
	\label{numbers}
	\end{gather}
where the notation $a \twoheadrightarrow b$ means that $a$ is divisible by $b$. Indeed, the equality $k_{F}(n)=k_{F}(n,0)$ follows from the additivity axiom for homology theories, the divisions $k_{F}(n,l)\,  |\,  k^{s}_{F}(n-l-1)$ hold since for an $l$-connected $CW$-complex $X$ the spectrum $\Sigma^{-(l+1)}\Sigma^{\infty}X$ is connective and $\mu_{F}$ commutes with suspension, and for $l\ge [\frac{n+1}{2}]$ the above divisions become equalities by the Freudenthal suspension theorem. One of the main questions of the paper is how one can find the \textit{minimal} positive integer $l$ such that 
$$k_{F}(n,l)=k^s_{F}(n-l-1)?$$ 

As is often the case in algebraic topology, it is more convenient to work not with arbitrary spectra but with $p$-local ones. So, let $F_{(p)}$ denote the localisation of $F$ at a prime number $p$. Then it is easy to check that $k_{F}(n,l)=\prod_{p\ \text{prime}}k_{F_{(p)}}(n,l)$ and $k^s_{F}(n)=\prod_{p\ \text{prime}}k^s_{F_{(p)}}(n)$. From now on let $p$ be a fixed prime number and $F$ be a fixed $p$-local spectrum of finite type (that is the homotopy groups of $F$ are finitely generated $\Z_{(p)}$-modules). One of the main results of the paper is the following theorem. 

\begin{Th}\label{Th1}
There exists a series of infinite loop spaces $\{\F^{n}_{(-r,0)}\}$ indexed by positive integers $n$ and $r$, distinguished cohomology classes $\mathfrak{f}^r_{n}\in H^{n}\bigl(\F^{n}_{(-r,0)}; \Hom_{\Z}(\pi_{r-1}F, \Q/\Z)\bigr)$ and a series of stable cohomology operations $\chi(\theta_r)\colon H^*\bigl(X; \Hom_{\Z}(\pi_{r-1}F, \Q/\Z)\bigr) \to H^{*+r}\bigl(X; \Hom_{\Z}(\pi_{0}F, \Q/\Z)\bigr)$ such that the following two conditions are equivalent:
\begin{enumerate}
\item $k_{F}(n, n-r-1)=k^s_{F}(r)$;
\item there exists $r'\le r$ such that $\chi(\theta_{r'})\bigl(\mathfrak{f}^{r'}_{n-r'}\bigr)\ne 0$ and  $k^s_{F}(r)=k^s_{F}(r')$. 
\end{enumerate}
\end{Th}

The theorem has surprising consequences. It turns out that if $r(n)$ denotes  the maximal positive integer such that $k_{F}\bigl(n, n-r(n)-1\bigr)=k^s_{F}\bigl(r(n)\bigr)$, then for $F=\MSO_{(p)},\MU_{(p)}$ or $\BU_{(p)}$ the ratio $\frac{n-r(n)-1}{n}$ tends to $0$ as $n$ tends to infinity. Also, the author's lower bound on $k_{\SO}(n)$, which we discussed above, can be obtained by a direct application of the theorem. For details and more examples see Section $3$.

Our proof of Theorem \ref{Th1} heavily relies on an interpretation of the numbers $k_{F}(n,l)$ and $k_{F}^s(n)$ in terms of Atiyah---Hirzebruch spectral sequences (AHSS for short). Namely, let $\{E^r_{s,t}(X); d^r_{s,t}\}$ denote the AHSS for the homology theory $F_*$ of a $CW$-complex or spectrum~$X$ and let $Z^r_{s,t}(X), B^r_{s,t}(X) \subset E^2_{s,t}(X)$ be the subgroups of cycles and boundaries correspondingly, so that $E^r_{s,t}(X)=\frac{Z^r_{s,t}(X)}{B^r_{s,t}(X)}$. Then the following holds. 

\begin{Prop}\label{ahss}
The number $k_{F}(n)$ (respectively, $k^{s}_{F}(n)$) is equal to the maximum of exponents of the groups $E^{2}_{n,0}(X)/Z^{\infty}_{n,0}(X)$ over all $CW$-complexes (respectively, connective spectra) $X$. 
\end{Prop}
\begin{proof}
Let  $\{\widetilde{E}^r_{s,t}(X); \widetilde{d}^r_{s,t}\}$ denote the AHSS for $\widetilde{F}_*$ of $X$. From Theorems $3.2$ and $3.3$ of \cite{Maunder} it follows that  $Z^{\infty}_{n,0}(X)\cong \widetilde{E}^{\infty}_{n,0}(X)$. Let $\mu'_{F}\colon \widetilde{F}_n(X) \to H_n\bigl(X;\pi_{0}(F)\bigr)$ denote the following composition of canonical homomorphisms:
\[
\xymatrix{
\widetilde{F}_n(X) \ar@{->>}[r]  & \widetilde{E}^{\infty}_{n,0}(X) \ar@{^{(}->}[r] & \widetilde{E}^2_{n,0}(X) \ar[r]^-{\cong} & H_n\bigl(X;\pi_{0}(F)\bigr).
}
\]
Then $\mu'_{F}$ is a natural transformation of the corresponding homology theories and for $X=\Sigma^{n}\mathbb{S}$ we have that $\mu'_{F}$ is the canonical isomorphism. So $\mu'_{F}\colon \widetilde{F} \to H\pi_{0}F$ induces the identity morphism on $\pi_{0}$ and thus coincides with~$\mu_{F}$. 
\end{proof}

Replacing $X$ with its factor by the appropriate skeleton, we obtain the following corollary. 

\begin{Cor}\label{Cor}
The number $k_{F}(n, n-r-1)$ (respectively, $k^{s}_{F}(r)$) is equal to the maximum of exponents of the groups $E^{2}_{n,0}(X; F)/Z^{r+1}_{n,0}(X; F)$ over all $CW$-complexes (respectively, spectra) $X$. 
\end{Cor}

Therefore the questions discussed above have a direct reformulation in terms of AHSSs. From this point of view, the core statement behind Theorem \ref{Th1} is the following. 

\begin{Th}\label{Th2}
For the same spaces $\{\F^{n}_{(-r,0)}\}$, classes $\mathfrak{f}^r_{n}$ and operations $\chi(\theta_r)$ as in Theorem \ref{Th1} the following two conditions are equivalent:
\begin{enumerate}
\item $\chi(\theta_r)\bigl(\mathfrak{f}^r_{n-r}\bigr)\ne 0$;
\item there exists a $CW$-complex $X$ and a homology class $x\in H_n(X;\pi_{0}F)$ such that $d^{r}_{n,0}\bigl(k^s_{F}(r-1)x\bigr)\ne 0$ in the AHSS for $F_*$ of $X$.
\end{enumerate}
\end{Th}

The paper is organised as follows. In Section $1$ we discuss the stable case and introduce operations~$\theta_r$. In Section $2$ we construct spaces $\F^{n}_{(-r,0)}$, classes $\mathfrak{f}^r_{n}$ and an operator $\chi$ and prove Theorems \ref{Th1} and \ref{Th2}. In Section~$3$ we discuss computational aspects of the theorems by establishing a connection between  this paper and the paper~\cite{R}. In Section~$4$ we give some concluding remarks.
 
I am grateful to my advisor Alexander Gaifullin for valuable discussions. 
 
\subsection{Stable case.} In this section we will describe one of the ways to compute the numbers $k^s_{F}(n)$. 

It turns out that in the stable category the numbers $k^s_{F}(n)$ can be interpreted in terms of a single cohomology class. Namely, let $\iota\in H^0(H\pi_{0}F, \pi_{0}F)$ be the fundamental class. Then the following holds. 

\begin{Lem}\label{L1}
The number $k^s_{F}(n)$ is the smallest positive integer such that the class $k^s_{F}(n)\iota$ serves as a cycle for all the differentials $d_r^{0,0}$ with $r\le n$ in the AHSS for $F^*$ of $H\pi_{0}F$. 
\end{Lem}

\begin{proof}
Let $A$ denote the group $\pi_{0}F$ and $A'$ be a finitely generated abelian group such that $A'\otimes \Z_{(p)}\cong A$. 
Also let $\iota'\in H^{0}(HA'; A)$ be the canonical class. Since the AHSSs for $F^*$ of $HA$ and of $HA'$ coincide (recall that the spectrum $F$ is $p$-local) it is sufficient to prove the lemma for the class $\iota'$ instead of $\iota$. 

Now let $X$ be a connective spectrum with a homology class $x\in H_n(X; A)$. Consider a finite cellular connective spectrum $X'$ of dimension $n$ with a homology class $x'\in H_n(X';A)$ and a map $f\colon X' \to X$ such that $f_*(x')=x$. Let $D$ denote the Spanier---Whitehead duality functor. It is easy to see that there exists a map $g \colon D(X')\to \Sigma^{-n}HA'$ such that $g^*(\Sigma^{-n}\iota')=D(x')$. By cellular approximation theorem the map~$g$ factorises through a map $\bar{g}\colon D(X') \to Sk^{0}\Sigma^{-n}HA'$ where $Sk^{0}$ stands for the $0$-th skeleton. As a result we obtain a map $D(\bar{g})\colon D(Sk^{0}\Sigma^{-n}HA') \to X'$ such that $D(\bar{g})_*D(\iota'_n)=x'$, where $\iota'_n$ is the pullback of $\Sigma^{-n}\iota'$ under the inclusion of the skeleton. 

Functoriality of AHSS and Proposition \ref{ahss} implies that $k^s_{F}(n)$ is the smallest positive integer such that $k^s_{F}(n)D(\iota'_n)$ is a permanent cycle in the AHSS for $F_*$ of $D(Sk^{0}\Sigma^{-n}HA')$. Since for a finite spectrum $X$ the AHSS for $F_*$ of $D(X)$ coincides with the AHSS for $F^*$ of $X$ with the reversed grading, we obtain that $k^s_{F}(n)$ is the smallest positive integer such that $k^s_{F}(n)\iota'_n$ is a permanent cycle in the AHSS for $F^*$ of $Sk^{0}\Sigma^{-n}HA'$. But this is already equivalent to the statement of the lemma. 
\end{proof}

Thereby the problem of computing the numbers $k^s_{F}(n)$ is closely related with the problem of understanding differentials $d_r^{0,0}$  in the AHSS for $F^*$ of $H\pi_{0}F$. Despite the fact that the latter problem is interesting and complicated, we will not discuss it here, but will assume that it has already been solved. Namely, we will need the following data:
\begin{enumerate}
\item the order of $\iota$ in the group $H^{0}\bigl(H\pi_{0}F, \pi_{0}F\bigr)/Z_r^{0,0}\bigl(H\pi_{0}F\bigr)$, which we will denote by $e_r$; 
\item an element $\theta_r\in H^{r}(H\pi_{0}F, \pi_{r-1}F)$ such that $\theta_r\in d_r^{0,0}(e_r\iota)$. 
\end{enumerate}

By Lemma \ref{L1} we have that $k^s_{F}(n)=e_{n+1}$, so the data $(1)$ is equivalent to the knowledge of the numbers~$k^{s}_{F}(n)$. The element $\theta_r$ provides us in a standard way with  the stable cohomological operation $\theta_r$, which acts on homology and cohomology groups: 
\[
\theta_r\colon H^n(X; \pi_{0}F) \to H^{n+r}(X; \pi_{r-1}F)\ \ \text{and} \ \ \theta_r\colon H_n(X; \pi_{0}F) \to H_{n-r}(X; \pi_{r-1}F) 
\]
where $X$ is a spectrum or  $CW$-complex. Taken together the data yields the following corollary. 

\begin{Cor}\label{Cor1}
For any homology class $x\in H_n(X;\pi_{0}F)$ of a $CW$-complex or connective spectrum $X$ in the AHSS for~$F_*$ of~$X$  we have that $e_rx\in Z^{r}_{n,0}$ and $\theta_r(x)\in d^{r}_{n,0}(e_rx)$. 
\end{Cor}
\begin{proof}
Since $X$ can be replaced by its suspension spectrum, the corollary needs only be proved for connective spectra. Using the same arguments as in the proof of Lemma \ref{L1}, it is easy to see that it is sufficient to prove the cohomological version of the statement for the class $\iota\in H^0(H\pi_{0}F, \pi_{0}F)$. But this is precisely our presumed data.  
\end{proof}

In the case $F=\BU_{(p)}$, Buchstaber \cite[Theorem 2.1]{Buch.2} proved that $e_r=p^{\left[\frac{r-2}{2(p-1)}\right]}$ and as the operations~$\theta_r$  one can take $\theta_r=\beta P^k$ for  $r=2k(p-1)+1$ and  $\theta_r=0$ for $r\ne1 \bmod 2(p-1)$.

\subsection{Proofs.} In this section we will prove Theorems \ref{Th1} and~\ref{Th2}. 

Let $\widehat{F}$ be a spectrum which represents the cohomology theory given by 
\[
\F^*(X)=\Hom_{\Z}\bigl(F_*(X), \Q/\Z \bigr)
\]
for any spectrum (and therefore $CW$-complex) $X$. Note that in particular $\pi_{-i}(\F)= \Hom_{\Z}(\pi_i(F), \Q/\Z)$. The spectrum $\widehat{F}$ is known as the Brown--Comenetz dual spectrum to~$F$ (see \cite{BC}). Now let 
\begin{itemize}
\item $\F_{(-r,0)}=\tau_{(-r,0)}\F$ be the spectrum $\F$ with killed homotopy groups in dimensions $s$ with $s\le-r$ and $s\ge 0$;
\item $\F^{n}_{(-r,0)}=\Omega^{\infty}\Sigma^{n+r-1}\F_{(-r,0)}$ be the $n$-th space of the $\Omega$-spectrum equivalent to $\Sigma^{r-1}\F_{(-r,0)}$;
\item $\mathfrak{f}^r_{n}\in H^{n}\bigl(\F^{n}_{(-r,0)};\pi_{-r+1}(\F)\bigr)$ be the canonical class. 
\end{itemize}

Observe that since the group $\Q/\Z$ is injective, the universal coefficient theorem implies that the following canonical map is an isomorphism:
\[
H^{n}\bigl(X; \Hom_{\Z}(A, \Q/\Z) \bigr) \xrightarrow{\cong} \Hom_{\Z}\bigl(H_{n}(X; A), \Q/\Z\bigr). 
\]
Now for a stable operation $\psi\colon H_{n}(X; A) \to H_{n-d}(X; B)$ of degree $d$  let 
\[
\chi(\psi)\colon H^{n-d}\bigl(X; \Hom_{\Z}(B, \Q/\Z)\bigr)\to H^{n}\bigl(X; \Hom_{\Z}(A, \Q/\Z)\bigr)
\]
be the operation obtained from $\psi$ by applying the functor $\Hom_{\Z}(-, \Q/\Z)$. Note that if $A\cong B\cong \Z/p\Z$, then $\psi$ is an element of the $\bmod \ p$ Steenrod algebra~$\mathcal{A}_{p}$ and the operator $\chi$ is the antipode in the Hopf algebra  $\mathcal{A}_{p}$ (this was proved at the beginning of  §4, Ch.3, \cite{Thom} and later reproved in \cite[Theorem 1.9(v)]{BC}). Also note that $\chi(\psi\circ\psi')=\chi(\psi')\circ\chi(\psi)$ when both sides are defined. From this it is not hard to express~$\chi(\psi)$ in standard terms in the case when $A$ and $B$ are finitely generated $\Z_{(p)}$ or $\Z$-modules. 

So, we have introduced all the notations which were used in the formulations of Theorems \ref{Th1} and \ref{Th2}. Now we need to do some preliminary work before going to proofs. 

Let $\{\E_r^{s,t}(X)=\frac{\Zz_r^{s,t}(X)}{\B_r^{s,t}(X)}; \dd_r^{s,t}\}$ denote the AHSS for~$\F^*$ of a $CW$-complex or spectrum $X$. Since \linebreak $\F^*(X)=\Hom_{\Z}\bigl(F_*(X), \Q/\Z \bigr)$ and $\Q/\Z$ is injective, we have the following description of it:
\[
\E_r^{s,t}(X)=\Hom_{\Z}\bigl(E^r_{s,t}(X), \Q/\Z\bigr); \ \dd_r^{s,t}=\Hom_{\Z}\bigl(d^r_{s,t}, \Q/\Z \bigr);
\]
\[
\Zz_r^{s,t}(X)=\Ann_{\Hom(E^2_{s,t}(X), \Q/\Z)}\bigl(B^r_{s,t}(X)\bigr); \ \B_r^{s,t}(X)=\Ann_{\Hom(E^2_{s,t}(X), \Q/\Z)}\bigl(Z^r_{s,t}(X)\bigr).
\]
The following lemma is a key ingredient for the proof of Theorem \ref{Th2}. 
\begin{Lem}\label{Lem2}
Let $X$ be a $CW$-complex or spectrum. If there exists a cohomology class $y\in H^{n-r}(X; \pi_{-r+1}\F)$ such that $y\in \Zz_r^{n-r, r-1}(X)$ and $\chi(\theta_r)(y)\ne 0$, then there exists a homology class $x\in H_n(X; \pi_{0}F)$ such that $d^r_{n,0}(e_rx)\ne~0$. Moreover, if all homology groups $H_{m}(X;\Z_{(p)})$ are finitely generated $\Z_{(p)}$-modules, then the converse statement also holds. 
\end{Lem}
\begin{proof}
Since 
\[
0\ne \chi(\theta_r)(y) \in H^n(X; \pi_{0}\F)\cong \Hom_{\Z}\bigl(H_n(X; \pi_{0}F), \Q/\Z  \bigr),
\]
then there exists $x\in H_n(X; \pi_{0}F)$ such that $\langle x, \chi(\theta_r)(y)\rangle \ne 0$. Since 
\[
y\in \Zz_r^{n-r, r-1}(X)=\Ann_{\Hom(E^2_{n-r,r-1}(X), \Q/\Z)}\bigl(B^r_{n-r,r-1}(X)\bigr),
\]
then for any representatives $x_1, x_2\in H_{n-r}(X; \pi_{r-1}F)$ of $d^r_{n,0}(e_rx)$, we have that $\langle x_1, y\rangle=\langle x_2, y\rangle$. By Corollary \ref{Cor1} we have that  $\theta_r(x)\in d^r_{n,0}(e_rx)$. So, since  $\langle \theta_r(x), y\rangle=\langle x, \chi(\theta_r)y\rangle\ne 0$, then  $0\not \in d^r_{n,0}(e_rx)$. 

Now let us prove the ``moreover'' part of the lemma. Since $F$ is $p$-local and of finite type, the group $H_{n-r}(X; \pi_{r-1}F)$ is finitely generated over $\Z_{(p)}$ and so has the form $T\oplus F$ where $T$ is a finite $p$-group and $F$ is a free $\Z_{(p)}$-module. Since differentials of AHSS are always of finite order we have that $B^r_{n-r, r-1}(X)\subset T$ and $d^r_{n,0}(e_rx)\in T/B^r_{n-r, r-1}(X)$. Since $T$ is finite it is easy to see that there exists a homomorphism \linebreak  $f\colon T \to \Q/\Z$ such that $f\bigl(B^r_{n-r, r-1}(X)\bigr)=0$ and $f(\tilde{x})\ne 0$ for all representatives $\tilde{x}\in T$ of $d^r_{n,0}(e_rx)$. Now consider a class $y\in H^{n-r}(X; \pi_{-r+1}\F)\cong\Hom_{\Z}(H_{n-r}(X; \pi_{r-1}F), \Q/\Z  \bigr)$ such that $y$ restricted to $T$ coincides with~$f$. Then $y\in \Zz_r^{n-r, r-1}(X)$ and for any representative $\tilde{x}\in H_{n-r}(X; \pi_{r-1}F)$ of $d^r_{n,0}(e_rx)$ we have that $\langle \tilde{x}, y\rangle\ne 0$. By Corollary \ref{Cor1} we know that $\theta_r(x)\in d^r_{n,0}(e_rx)$. Therefore $0\ne \langle \theta_r(x), y\rangle=\langle x, \chi(\theta_r)y\rangle$ and so $\chi(\theta_r)y\ne 0$. 
\end{proof}

In \cite{Maunder}, C.~R.~F. Maunder characterized the pages and differentials of the AHSS for a cohomology theory $Y^*$ using  the Postnikov system of the spectrum $Y$ (see especially Theorem $3.2$, Theorem~$3.4$ and their proofs). In particular, his results show that the class $\mathfrak{f}^r_{n}\in H^{n}(\F^{n}_{(-r,0)};\pi_{-r+1}\F)$ lies in the subgroup $\Zz_r^{n, r-1}(\F^{n}_{(-r,0)})$. Moreover, for any cohomology class $y\in H^{n}(X; \pi_{-r+1}\F)$ of a $CW$-complex $X$ such that $y\in\Zz_r^{n, r-1}(X)$ there exists a map $f\colon X \to \F^{n}_{(-r,0)}$ satisfying $f^*(\mathfrak{f}^r_{n})=y$. 

\begin{proof}[Proof of Theorem \ref{Th2}.]
As was mentioned above, from \cite{Maunder} it follows that $\mathfrak{f}^r_{n-r} \in \Zz_r^{n-r, r-1}(\F^{n-r}_{(-r,0)})$. Since $k^s_{F}(r-1)=e_r$, the implication $(1)\to (2)$ completely follows from Lemma \ref{Lem2}. 

Now let us prove that $(2)$ implies $(1)$. Let $x\in H_n(X;\pi_{0}F)$ be a homology class of a CW-complex $X$ such that $d^{r}_{n,0}(e_{r}x)\ne 0$. Since we can realize $x$ by a homology class of a finite $CW$-complex, functoriality of AHSS implies that we can assume that $X$ is finite. Now by Lemma \ref{Lem2} there exists a cohomology class $y\in H^{n-r}(X; \pi_{-r+1}\F)$ such that $y\in \Zz_r^{n-r, r-1}(X)$ and $\chi(\theta_r)(y)\ne 0$. As was mentioned above, from \cite{Maunder} it follows that there exists a map $f\colon X \to \F^{n-r}_{(-r,0)}$ such that $f^*(\mathfrak{f}^r_{n-r})=y$. Now functoriality with respect to~$f$ implies that $\chi(\theta_r)(\mathfrak{f}^r_{n-r})\ne 0$ since $\chi(\theta_r)(y)\ne 0$. So the theorem is proved. 
\end{proof}

Note that in the proof we did not use the fact that the spectrum $F$ is $p$-local. So the theorem remains to be true if $F$ is an arbitrary spectrum of finite type. 

For the proof of Theorem \ref{Th1} we need the following auxiliary statement.

\begin{Prop}\label{Prop}
For any number $r$ we have that either $e_{r+1}=e_r$ or $e_{r+1}=pe_{r}$. 
\end{Prop}
\begin{proof}
Recall that if $A, B$ are finitely generated $p$-local abelian groups and $\psi \in H^{m}(HA; B)$ is a non-zero cohomology class of degree $m\ge2$, then $p\psi=0$ (see \cite[Th.~2]{Cartan}). Now recall that $e_r$ is the order of $\iota$ in the group $H^{0}\bigl(H\pi_{0}F, \pi_{0}F\bigr)/Z_r^{0,0}\bigl(H\pi_{0}F\bigr)$ and $\theta_r\in H^{r}\bigl(H\pi_{0}F, \pi_{r-1}F\bigr)$ is such that $\theta_r\in d_r^{0,0}(e_r\iota)$. If $d_r^{0,0}(e_r\iota)=0$, then clearly $e_{r+1}=e_r$. Now if $d_r^{0,0}(e_r\iota)\ne0$, then $0=p\theta_r\in d_r^{0,0}(pe_r\iota)$ and so $e_{r+1}\le pe_{r}$, but on the other hand the ideal of non-invertible elements in~$\Z_{(p)}$ is generated by $p$ and so $e_{r+1}\ge pe_{r}$. Therefore $e_{r+1}=pe_{r}$.
\end{proof}

Note that the same argument was used by Buchstaber in the proof of \cite[Th.~3]{Buch.1} and in \cite[Rem.~2]{Buch.1}. Now we are ready to give a proof of Theorem \ref{Th1}. 

\begin{proof}[Proof of Theorem \ref{Th1}.] 
By Theorem \ref{Th2}, the condition $\chi(\theta_{r'})(\mathfrak{f}^{r'}_{n-r'})\ne 0$ is equivalent to the existence of a homology class $x\in H_{n}(X;\pi_{0}F)$ of a $CW$-complex $X$ with the property that $d_{r'}^{n,0}(e_{r'}x)\ne 0$. Therefore the order of $x$ in the group 
$
E^{2}_{n,0}(X)/Z^{r'+1}_{n,0}(X)
$
 is at least $pe_{r'}$. On the other hand, by Corollary \ref{Cor1} it is at most $e_{r'+1}$. So, Proposition \ref{Prop} implies that $pe_{r'}=e_{r'+1}$ and thus the order of $x$ in $E^{2}_{n,0}(X)/Z^{r'+1}_{n,0}(X)$ is equal to~$k^{s}_{F}(r')$. Therefore Corollary \ref{Cor} implies that the condition $\chi(\theta_{r'})(\mathfrak{f}^{r'}_{n-r'})\ne 0$ is equivalent to the condition that $k_{F}(n, n-r'-1)=k^{s}_{F}(r')$ and $k^{s}_{F}(r')>k^{s}_{F}(r'-1)$. Finally, Diagram \eqref{numbers} implies that the latter conditions are equivalent to the condition that for any $r$ such that $k^{s}_{F}(r')=k^{s}_{F}(r)$ one has that $k_{F}(n, n-r-1)=k^{s}_{F}(r)$. 
\end{proof}

\subsection{Computational aspects and a connection with \cite{R}.} Let us take for $F$ the spectrum $\BU_{(p)}$ of complex $K$-theory, localised at a fixed prime number $p$. Then we have the following theorem, which was essentially proved in \cite{R}. 

\begin{Th}\label{Th3}
For $F=\BU_{(p)}$, if  
\[
r(n)=2\Bigr(\sum_{i=1}^{\infty} \left[\frac{n-1}{2p^i}\right]\Bigl)(p-1)+1
\]
then $\chi(\theta_{r(n)})(\mathfrak{f}^{r(n)}_{n-r(n)})\ne 0$. 
\end{Th}
\begin{proof}
Recall that the spectrum $\BU$ is Anderson self-dual (see \cite{A.duality}) and so we have the exact sequence $\BU_{(p)} \to \BU_{(0)} \to \widehat{\BU_{(p)}}$, where $\BU_{(0)}$ stands for the rationalization of $\BU$. Using this and the Bott periodicity theorem it is not hard to check that the space $R_{p}(q,s)$ introduced in \cite[Sect. 4]{R} is homotopy equivalent to the space $\widehat{F}^{q-1}_{(-2s-1, 0)}$ constructed above. In particular, the class $\tilde{u}_{q,s}$ introduced in \cite{R} is equal to the class $\beta\mathfrak{f}^{2s+1}_{q}$, where $\beta$ is the Bockstein homomorphism for the short exact sequence $0\to \Z_{(p)} \to \Q \to \Q/\Z_{(p)}\to 0$. Now the statement completely follows from Lemma 10(3) and Lemma 7 in \cite{R}. 
\end{proof}

\begin{Cor}
The statement of Theorem \ref{Th3} also holds for $F=\MSO_{(p)}$ with $p$ an odd prime as well as for $F=\MU_{(p)}$ with any prime $p$.  
\end{Cor}
\begin{proof}
Let us first consider the case $F=\MU_{(p)}$. Let $\mu_{c}\colon \MU_{(p)} \to \BU_{(p)}$ be the Conner---Floyd map and recall that  $\mu_{c*}\colon \pi_{0}\MU_{(p)} \to \pi_{0}\BU_{(p)}$ is an isomorphism. Also, by \cite{Buch.2} we have that $k^s_{\BU}(n)=k^s_{\MU}(n)$ for all $n$. Now the statement completely follows from Theorem \ref{Th2} by the functoriality of AHSS along  $\mu_{c*}$. 

Now consider $F=\MSO_{(p)}$, where $p$ is odd. From Theorem 1.5 of \cite{BP} it follows that for any $n$ we have that $k^s_{\MSO_{(p)}}(n)=k^s_{\BP}(n)=k^s_{\MU_{(p)}}(n)$, where $\BP$ stands for the Brown---Peterson spectrum at the prime $p$. Also Theorem 1.3 of \cite{BP}  implies that $\BP$ is a retract of $\MSO_{(p)}$ as well as of $\MU_{(p)}$ and these retractions induce isomorphisms on $\pi_{0}$. Now the statement completely follows from Theorem \ref{Th2} using the functoriality of AHSS along these retraction maps.
\end{proof}

Now let $F$ be either $\BU_{(p)}$ or $\MU_{(p)}$ for any prime $p$, or $\MSO_{(p)}$ for an odd prime $p$. In \cite{Buch.2} it was proved that $k^{s}_{F}(n)=p^{\left[\frac{n-1}{2(p-1)}\right]}$. Therefore, the above results together with Theorem \ref{Th1} imply that 
\[
k_{F}^{s}(r(n))=k_{F}(n, n-r(n)-1)=p^{\sum_{i=1}^{\infty} \left[\frac{n-1}{2p^i}\right]}.
\]
It is easy to check (see also the proof of \cite[Cor. 2]{R}) that 
\[
\lim_{n\to \infty}\frac{n-r(n)-1}{n}=0; \ \ \  \lim_{n\to \infty} \frac{\log_{p}\bigl(k_{F}(n)\bigr)}{\log_{p}\bigl(k_{F}^{s}(r(n))\bigr)}= \lim_{n\to \infty}\frac{\log_{p}\bigl(k_{F}(n)\bigr)}{\log_{p}\bigl(k_{F}^{s}(n)\bigr)}=1. 
\]
Also, in \cite[Th. 12]{R} it was proved that for the above spectra $k_{F}(n)=p^{\sum_{i=1}^{\infty} \left[\frac{n-1}{2p^i}\right]}$ for $n<2p^2+2p$. So Theorem \ref{Th1} indeed yields quite strong lower bounds on $k_{F}(n)$. 

Let us mention that in order to obtain some information about $k_{F}(n)$ for a spectrum $F$ it is not necessary to apply Theorem \ref{Th1} directly to $F$. Namely, if $f\colon F\to F'$ is a map of spectra such that \linebreak $f_{*}\colon \pi_{0}F \to \pi_{0}F'$ is an isomorphism, then one obviously has that $k_{F}(n)$ is divisible by $k_{F'}(n)$ and $k^{s}_{F}(n)$ is divisible by $k^{s}_{F'}(n)$ for any $n$. Now one can use Theorem \ref{Th1} for $F'$ in order to obtain a lower bound on $k_{F'}(n)$ and therefore on $k_{F}(n)$. Applying this to the map $\mathbb{S}_{(p)} \to \MU_{(p)}$ we obtain that  $k_{\mathbb{S}_{(p)}}(n)$ is divisible by $p^{\sum_{i=1}^{\infty} \left[\frac{n-1}{2p^i}\right]}$. As was shown in \cite[Cor.3]{R} the ratio $\frac{\log_{p}\left(k_{\mathbb{S}_{(p)}}(n)\right)}{\sum_{i=1}^{\infty} \left[\frac{n-1}{2p^i}\right]}$ also tends to $1$ as~$n$ tends to $\infty$. On the other hand, we do not know how to find a sufficiently small $l$ such that $k_{\mathbb{S}_{(p)}}(n,l)=k^{s}_{\mathbb{S}_{(p)}}(n-l-1)$.  

\subsection{Conclusion. } In attempting to compute the numbers $k_{F}(n)$, one usually faces the following two difficulties. The first is the absence of an analog of Lemma \ref{L1} in the unstable category (due to the fact that the Alexander duality is not an anti-equivalence on the category of finite $CW$-complexes). Without a universal (co)homology class on which all the computations can be done, it is unclear how to obtain precise information about $k_{F}(n)$. The second is that a direct computation of differentials in a homological AHSS in the unstable category seems to be a hopeless problem. 

Theorem \ref{Th1} allows us to work around both of these difficulties simultaneously. We work with fixed cohomology classes and we need to compute the action of ordinary (but not higher kind) cohomological operations. Even though we cannot find the value of $k_{F}(n)$ by this method, we can understand how far the computations in the stable category can be transferred to the unstable setting. 

The most significant idea in the proofs of Theorem \ref{Th1} and Theorem \ref{Th2} is the simultaneous use of the Spanier---Whitehead duality, in order to obtain Corollary \ref{Cor1} from computations on a fixed cohomology class $\iota$, and the Brown---Comenetz duality, in order to prove Lemma \ref{Lem2} which, presumably, is useful already in the stable category but especially useful in the unstable one. It seems that this idea is new and we hope that it can be used to make progress in other semi-stable problems in topology.

\end{document}